47A48    Operator colligations (= nodes), vessels, 
linear systems, characteristic functions, realizations, etc.

\documentclass{article}

\usepackage{amssymb,amsmath,theorem,euscript}

\newcounter{sec}

\newcounter{punct}[sec]

\def\punct{\refstepcounter{punct}{\arabic{sec}.\arabic{punct}.  }}

\def\COUNTERS{\addtocounter{sec}{1}
              \setcounter{punct}{0}
          \setcounter{equation}{0}
          \setcounter{theorem}{0}
                  }

\newtheorem{theorem}{Theorem}[sec]
\newtheorem{proposition}[theorem]{Proposition}

\newtheorem{lemma}[theorem]{Lemma}

\newtheorem{corollary}[theorem]{Corollary}

 \def\ov{\overline}

 \newcommand{\rk}{\mathop {\mathrm {rk}}\nolimits}
\renewcommand{\Im}{\mathop {\mathrm {Im}}\nolimits}
\renewcommand{\Re}{\mathop {\mathrm {Re}}\nolimits}

\begin{document}

\def\OO{\mathrm{O}}
\def\GLO{\mathrm{GLO}}
\def\Coll{\mathrm{Coll}}
\def\kappa{\varkappa}

\def\R{\mathbb{R}}
\def\C{\mathbb{C}}

\def\la{\langle}
\def\ra{\rangle}

 \def\cA{\mathcal A}
\def\cB{\mathcal B}
\def\cC{\mathcal C}
\def\cD{\mathcal D}
\def\cE{\mathcal E}
\def\cF{\mathcal F}
\def\cG{\mathcal G}
\def\cH{\mathcal H}
\def\cJ{\mathcal J}
\def\cI{\mathcal I}
\def\cK{\mathcal K}
 \def\cL{\mathcal L}
\def\cM{\mathcal M}
\def\cN{\mathcal N}
 \def\cO{\mathcal O}
\def\cP{\mathcal P}
\def\cQ{\mathcal Q}
\def\cR{\mathcal R}
\def\cS{\mathcal S}
\def\cT{\mathcal T}
\def\cU{\mathcal U}
\def\cV{\mathcal V}
 \def\cW{\mathcal W}
\def\cX{\mathcal X}
 \def\cY{\mathcal Y}
 \def\cZ{\mathcal Z}
\def\0{{\ov 0}}
 \def\1{{\ov 1}}
 \def\frA{\mathfrak A}
 \def\frB{\mathfrak B}
\def\frC{\mathfrak C}
\def\frD{\mathfrak D}
\def\frE{\mathfrak E}
\def\frF{\mathfrak F}
\def\frG{\mathfrak G}
\def\frH{\mathfrak H}
\def\frI{\mathfrak I}
 \def\frJ{\mathfrak J}
 \def\frK{\mathfrak K}
 \def\frL{\mathfrak L}
\def\frM{\mathfrak M}
 \def\frN{\mathfrak N} \def\frO{\mathfrak O} \def\frP{\mathfrak P} \def\frQ{\mathfrak Q} \def\frR{\mathfrak R}
 \def\frS{\mathfrak S} \def\frT{\mathfrak T} \def\frU{\mathfrak U} \def\frV{\mathfrak V} \def\frW{\mathfrak W}
 \def\frX{\mathfrak X} \def\frY{\mathfrak Y} \def\frZ{\mathfrak Z} \def\fra{\mathfrak a} \def\frb{\mathfrak b}
 \def\frc{\mathfrak c} \def\frd{\mathfrak d} \def\fre{\mathfrak e} \def\frf{\mathfrak f} \def\frg{\mathfrak g}
 \def\frh{\mathfrak h} \def\fri{\mathfrak i} \def\frj{\mathfrak j} \def\frk{\mathfrak k} \def\frl{\mathfrak l}
 \def\frm{\mathfrak m} \def\frn{\mathfrak n} \def\fro{\mathfrak o} \def\frp{\mathfrak p} \def\frq{\mathfrak q}
 \def\frr{\mathfrak r} \def\frs{\mathfrak s} \def\frt{\mathfrak t} \def\fru{\mathfrak u} \def\frv{\mathfrak v}
 \def\frw{\mathfrak w} \def\frx{\mathfrak x} \def\fry{\mathfrak y} \def\frz{\mathfrak z} \def\frsp{\mathfrak{sp}}
 \def\bfa{\mathbf a} \def\bfb{\mathbf b} \def\bfc{\mathbf c} \def\bfd{\mathbf d} \def\bfe{\mathbf e} \def\bff{\mathbf f}
 \def\bfg{\mathbf g} \def\bfh{\mathbf h} \def\bfi{\mathbf i} \def\bfj{\mathbf j} \def\bfk{\mathbf k} \def\bfl{\mathbf l}
 \def\bfm{\mathbf m} \def\bfn{\mathbf n} \def\bfo{\mathbf o} \def\bfp{\mathbf p} \def\bfq{\mathbf q} \def\bfr{\mathbf r}
 \def\bfs{\mathbf s} \def\bft{\mathbf t} \def\bfu{\mathbf u} \def\bfv{\mathbf v} \def\bfw{\mathbf w} \def\bfx{\mathbf x}
 \def\bfy{\mathbf y} \def\bfz{\mathbf z} \def\bfA{\mathbf A} \def\bfB{\mathbf B} \def\bfC{\mathbf C} \def\bfD{\mathbf D}
 \def\bfE{\mathbf E} \def\bfF{\mathbf F} \def\bfG{\mathbf G} \def\bfH{\mathbf H} \def\bfI{\mathbf I} \def\bfJ{\mathbf J}
 \def\bfK{\mathbf K} \def\bfL{\mathbf L} \def\bfM{\mathbf M} \def\bfN{\mathbf N} \def\bfO{\mathbf O} \def\bfP{\mathbf P}
 \def\bfQ{\mathbf Q} \def\bfR{\mathbf R} \def\bfS{\mathbf S} \def\bfT{\mathbf T} \def\bfU{\mathbf U} \def\bfV{\mathbf V}
 \def\bfW{\mathbf W} \def\bfX{\mathbf X} \def\bfY{\mathbf Y} \def\bfZ{\mathbf Z} \def\bfw{\mathbf w}
 \def\R {{\mathbb R }} \def\C {{\mathbb C }} \def\Z{{\mathbb Z}} \def\H{{\mathbb H}} \def\K{{\mathbb K}}
 \def\N{{\mathbb N}} \def\Q{{\mathbb Q}} \def\A{{\mathbb A}} \def\T{\mathbb T} \def\P{\mathbb P} \def\G{\mathbb G}
 \def\bbA{\mathbb A} \def\bbB{\mathbb B} \def\bbD{\mathbb D} \def\bbE{\mathbb E} \def\bbF{\mathbb F} \def\bbG{\mathbb G}
 \def\bbI{\mathbb I} \def\bbJ{\mathbb J} \def\bbL{\mathbb L} \def\bbM{\mathbb M} \def\bbN{\mathbb N} \def\bbO{\mathbb O}
 \def\bbP{\mathbb P} \def\bbQ{\mathbb Q} \def\bbS{\mathbb S} \def\bbT{\mathbb T} \def\bbU{\mathbb U} \def\bbV{\mathbb V}
 \def\bbW{\mathbb W} \def\bbX{\mathbb X} \def\bbY{\mathbb Y} \def\kappa{\varkappa} \def\epsilon{\varepsilon}
 \def\phi{\varphi} \def\le{\leqslant} \def\ge{\geqslant}

\def\Gms{\mathrm{Gms}}
\def\Pol{\mathrm{Pol}}
\def\O{\mathrm O}

\def\zigzag{\rightsquigarrow}

\def\sm{\smallskip}

\begin{center}
\Large\bf
Symmetries of Gaussian measures and operator colligations 

\bigskip

\large\sc
Yury Neretin%
\footnote{Supported by grants FWF, P22122, P19064,  and by RosAtom, contract H.4e.45.90.11.1059.}
\end{center}

{\small Consider an infinite-dimensional linear space equipped with a Gaussian measure and the group $\GLO(\infty)$
 of linear transformations that send the measure to equivalent one. Limit points 
of $\GLO(\infty)$ can be regarded as 'spreading' maps (polymorphisms). We show that the closure
of $\GLO(\infty)$ in the semigroup of polymorphisms contains a certain semigroup of operator colligations  
and write explicit formulas for action of operator colligations by polymorphisms
of the space with Gaussian measure.}

\section{Introduction. Polymorphisms,  Gaussian measures, and colligations}

\COUNTERS

{\bf \punct The group $\Gms(M)$.}
Let $M=(M,\mu)$ be a Lebesgue space $M$ with a probability measure $\mu$ (\cite{Roh}, see, also
\cite{MI}),
let $L^p(M,\mu)$ be the space of measurable functions on $M$ with norm
$$
\|f\|_p=\Bigl(\int_M |f(m)|^p\,d\mu(m)\Bigr)^{\frac 1p}
,\qquad
\text{where $1\le p\le \infty$}.
$$

Denote by $\Gms(M)$ the group of all bijective a.s. maps $M\to M$ that send the measure $\mu$
 to an equivalent measure.
For  $g\in \Gms(M)$ we denote by $g'(m)$ the Radon--Nikodym derivative   of $g$.

\smallskip

Fix $\lambda\in\C$ lying in the strip $0\le \Re\lambda\le 1$,
\begin{equation}
\lambda=\tfrac 1p+is,\qquad \text {where $1\le p\le \infty$, $s\in\R$}
\label{eq:lambda}
.
\end{equation}
For any $g\in\Gms(M)$ we define the linear operator $T_\lambda(g)$ by
\begin{equation}
T_\lambda(g)f(m)=f(mg)g'(m)^\lambda
\label{eq:T-lambda}
.
\end{equation}
Evidently, the operators $T_\lambda(g)$ form a representation of the group $\Gms(M)$ by
isometric operators in the Banach space $L^p(M,\mu)$. For $p=2$ we get a unitary representation
in $L^2(M,\mu)$.

Polymorphims, which are introduced below, are "limit points" of the group $\Gms(M)$.


\smallskip

{\bf\punct  Gaussian measures.} Consider $\R$ equipped with the Gaussian measure
$\frac 1{\sqrt{2\pi}}e^{-x^2/2}\,dx$. 
Let $n=1$, 2, \dots, $\infty$. Denote by $\R^\omega$ the product  of 
$n$ copies of $\R$ equipped with
the product measure $\mu_\omega=\mu\times\mu\times\dots$. 
We denote elements of $\R^\omega$ by $x=\begin{pmatrix}x_1,x_2,\dots\end{pmatrix}$.

\begin{proposition}
\label{pr:linear}
If $\sum b_j^2<\infty$, then the series $\sum b_j x_j$ converges a.s. on $\R^\infty$ with respect to
the measure $\mu_\infty$.
\end{proposition}

This is a special case of the Kolmogorov--Hinchin theorem about series of independent random 
variables, see, e.g., \cite{Shiryaev}.

\smallskip


{\bf\punct Groups of symmetries of Gaussian measures.}
Denote by $\OO(\infty)$ the infinite-dimensional orthogonal group, i.e.,
 the group of all infinite real matrices $A$ satisfying the conditions
$$
AA^t=A^tA=1
,$$
where $^t$ denotes the transposition.

For an invertible real infinite matrix $A$ we consider the
polar decomposition $A=S U$, where $U\in\OO(\infty)$, and $S$ is a positive self-adjoint operator.
We define the
%
%
%
group $\GLO(\infty)$ consisting of matrices $A=SU$ such that $S-1$ 
is a Hilbert--Schmidt%
\footnote{An operator $T$ is Hilbert--Schmidt, if $\sum_{ij}|t_{ij}|^2<\infty$, see, e.g., \cite{RS1}}
 operator. Equivalently, we can represent
$A$ as $A=\exp(T)U$, where $U\in\OO(\infty)$ and $T$ is a Hilbert--Schmidt self-adjoint operator.

\smallskip

Thus the set $\GLO(\infty)$ is the product of $\OO(\infty)$ and the 
space of self-adjoint Hilbert--Schmidt matrices. We take the
 weak operator topology%
\footnote{See e.g., \cite{RS1}.} on $\OO(\infty)$
 and the natural topology on the space of
 Hilbert--Schmidt matrices%
\footnote{See e.g. \cite{RS1}.} . We equip $\GLO(\infty)$
with the topology of product. Then $\GLO(\infty)$ is a topological group with respect
to this topology
({\it the Shale topology,} \cite{Sha}).

\smallskip

Consider an infinite matrix $A=\{a_{ij}\}$. Apply it to a vector $x\in\R^\infty$,
\begin{equation}
xA=\begin{pmatrix}x_1&x_2&\dots\end{pmatrix}
\begin{pmatrix}a_{11}&a_{12}&\dots\\
a_{21}&a_{22}&\dots\\
\vdots&\vdots&\ddots
\end{pmatrix}
=
\begin{pmatrix}
\sum x_i a_{i1} & \sum x_ia_{i2}& \dots
\end{pmatrix}
\label{eq:xA}
\end{equation}

Let $A$ be an operator bounded in the space $\ell_2$.
By Proposition \ref{pr:linear} 
the vector $xA$  is defined for almost all $x\in(\R^\infty,\mu_\infty)$.

\begin{theorem}
{\rm a)} For $A\in\OO(\infty)$ the map $x\mapsto xA$ preserves measure $\mu_\infty$.

\smallskip

{\rm b)} For $A\in\GLO(\infty)$, the map $x\mapsto xA$ is defined a.s. on $(\R^\infty,\mu_\infty)$ 
and
sends the measure $\mu_\infty$ to an equivalent measure
$\mu(xA)$.

\smallskip

{\rm c)} Let $A=(1+T)U$, where $A\in \OO(\infty)$ and $T$ is in the trace
class%
\footnote{See, \cite{RS1}.}. Then the Radon--Nikodym derivative is given by the formula
\begin{multline}
\frac{d\mu(xA)}{d\mu(x)}=
|\det A|\cdot \exp(-\tfrac12 \la xA,xA\ra+\tfrac12\la x,x\ra)
:=\\:=
|\det(1+T)|\cdot \exp\bigr(-\la xT,x\ra-\tfrac12\la xT,xT\ra\bigl)
\label{eq:RN}
\end{multline}

\smallskip

{\rm d)} Let $A=1+T$, where $T$ is a diagonal matrix with entries 
$t_j>-1$ satisfying $\sum_j  t_j^2<\infty$. Then the Radon--Nikodym derivative is given by
$$
\prod_{j=1}^\infty (1+t_j) e^{-(2t_j+t_j^2)x_j^2/2}
,$$
the product converges a.s. on $(\R^\infty,\mu_\infty)$.

\smallskip

{\rm e)} For $A$, $B\in\GLO(\infty)$ the identity
$$
(xA)B=x(AB)
$$
holds a.s. on $(\R^\infty,\mu)$.
\end{theorem}

The theorem is a reformulation of the Feldman--Hajeck Theorem on equivalence of Gaussian measures
(see, e.g., \cite{Kuo}, \cite{Bog}), the most comprehensive  exposition is in \cite{Shilov}.

\sm

{\sc Remark.}
For $A\in\GLO_1(\infty)$, the absolute value of determinant
 $|\det(A)|:=|\det(1+T)|$ is well-defined 
(see, e.g, \cite{Ner-book}),
it satisfies 
$$
|\det (A_1A_2)|=|\det(A_1)|\cdot |\det(A_2|
.$$
The $\det(A)$ makes no sence. \hfill $\square$

\sm

{\sc Remark.} In our definition the action is defined a.s, and the identity $x(AB)=(xA)B$
also is valid a.s. The removing of "a.s." is impossible, the group $\O(\infty)$ can not act 
pointwise by  measure preserving transformations, see \cite{Gla}.
 \hfill $\square$

\smallskip



{\bf\punct Polymorphisms (spreading maps)}, for  details, see \cite{Ner-poli}.
 \cite{Ner-book}, \cite{Ner-match}).
  Denote by $\R^\times$ the multiplicative group of positive real numbers,
denote by $t$ the coordinate on $\R^\times$, by $\alpha*\beta$ we denote the convolution
of measures on $\R^\times$.
Let $M=(M,\mu)$, $N=(N,\nu)$ be  Lebesgue spaces with  probability measures.
 A {\it polymorphism}%
\footnote{These objects were introduced in \cite{Ner-bist}, see also
\cite{Ner-book}. The term was proposed be Vershik \cite{Ver}, who used it for measures on $M\times N$,
see also "bistochastic kernels" from \cite{Kre}. On some appearances of polymorphisms in variation
problems and mathematical hydrodynamics, see \cite{Bre}.}
$\frP:(M,\mu)\zigzag (N,\nu)$  
is a measure $\frP=\frP(m,n,t)$ on $M\times N\times\R^\times$ satisfying two conditions:

\smallskip

a) the projection of $\frP(m,n,t)$ to $M$ is $\mu$;

\smallskip

b) the projection of $t\cdot\frP(m,n,t)$ to $N$ is $\nu$. 

\smallskip

We denote by $\Pol(M,N)$ the set of all polymorphisms $(M,\mu)\zigzag (N,\nu)$.

\smallskip

There is a well-defined associative  multiplication
$$
\Pol(M,N)\times \Pol(N,K)\to \Pol(M,K)
$$


{\bf\punct Convergence of polymorphisms.}
For $\frP\in\Pol(M,N)$ and measurable subsets $A\subset M$, $B\subset N$ we 
consider the projection $A\times B\times\R^\times\to\R^\times$
 and denote by $\frp[A\times B]$ 
the pushforward of $\frP$ under
this projection. 

We say that a sequence $\frP_j\in\Pol(M,N)$ {\it converges} 
to $\frP$ if for any $A\subset M$, $B\subset N$ we 
have weak convergences
$$
\frp[A\times B]\to\frp[A,\times B],\qquad
t\cdot\frp_j[A\times B]\to t\cdot \frp[A\times B]
.$$ 

\begin{proposition}
The product of polymorphisms is separately continuous,
i.e. if $\frP_j$ converges to $\frP$ in $\Pol(M,N)$ and $\frQ_j$
 converges to $\frQ$ in $\Pol(N,K)$,
then $\frQ\diamond \frP_j$ converges to $\frQ\diamond \frP$ and 
$\frQ_j\diamond \frP$ converges to $\frQ\diamond \frP$.
\end{proposition}

Note that there is no joint continuity, generally $\frQ_j\frP_j$ 
does not converge to $\frQ\diamond \frP$.

\sm


{\bf\punct Embedding $\frI:\Gms(M)\to\Pol(M,M)$.}
Now let a measure $\mu$ on $M$ be continuous. We consider the embedding 
\begin{equation}
\frI:\Gms(M)\to\Pol(M,M)
\label{eq:Gms-Pol}
\end{equation} 
given by the following way.
Take the map $M\mapsto M\times M\times\R^\times$ 
given by $m\mapsto \bigl(m,g(m),g'(m)\bigr)$. Then the
pushforward of the measure $\mu$ is a polymorphism 
$\frI(g):M\to M$.

\begin{proposition} {\rm(\cite{Ner-bist}, \cite{Ner-poli})}
The group $\Gms(M)$ is dense in $\Pol(M,M)$.
\end{proposition}


\smallskip



{\bf\punct  Formulation of problem.} We wish to describe the closure
of $\GLO(\infty)$ in the semigroup of polymorphisms%
\footnote{The closure of $\OO(\infty)$ gives action of the semigroup
of all contractive  linear operators by polymorphisms of $\R^\infty$, see  Nelson \cite{Nel}, .}
 of $\R^\infty$.  Our solution is not final, we show a large semigroup
(see the next subsection) in this closure.


\smallskip

{\bf\punct Operator colligations.}
 Fix $\omega=0$, 1,  \dots, $\infty$. Denote by $\GLO(\omega+\infty)$
the group consisting of $(\omega+\infty)\times(\omega+\infty)$ matrices $g$
 that are elements of the group
$\GLO$  (i.e, $\GLO(\omega+\infty)$ is another notation for $\GLO(\infty)$).
 Consider the subgroup $\OO(\infty)\subset \GLO(\omega+\infty)$
consisting of block $(\omega+\infty)\times(\omega+\infty)$
 matrices $\begin{pmatrix}1&0\\0&u \end{pmatrix}$, 
where $u$ is an orthogonal matrix.

We say that
an {\it operator colligation} is an 
element $g$ of $\GLO(\omega+\infty)$ defined up to the equivalence
$$
g\sim h_1 g h_2,\qquad \text{where $h_1$, $h_2\in\OO(\infty)$},
$$
or, in more details,
$$
\begin{pmatrix}\alpha&\beta\\ \gamma&\delta\end{pmatrix}
\sim
\begin{pmatrix}1&0\\0&u \end{pmatrix}
\begin{pmatrix}\alpha&\beta\\ \gamma&\delta\end{pmatrix}
\begin{pmatrix}1&0\\0&v \end{pmatrix}
$$ 
where $u$, $v$ are orthogonal matrices. Denote by $\Coll(\omega)$ the set 
of all operator colligations. In other words, $\Coll(\omega)$
 is the double coset space
$$
\Coll(\omega)=
\OO(\infty)\setminus \GLO(\omega+\infty) /\OO(\infty)
.$$

The {\it product of operator colligations} is defined by the formula
$$
\begin{pmatrix}\alpha&\beta\\ \gamma&\delta\end{pmatrix}
\circ
\begin{pmatrix}\phi&\psi\\\theta&\kappa \end{pmatrix}:=
\begin{pmatrix}\alpha&\beta&0\\ \gamma&\delta&0\\
0&0&1
\end{pmatrix}
\begin{pmatrix}\phi&0&\psi\\ 0&1&0\\\theta&0&\kappa \end{pmatrix}=
\begin{pmatrix}
\alpha\phi&\beta&\alpha\psi\\
\gamma\phi&\delta&\gamma\psi\\
\theta&0&\kappa
\end{pmatrix}
$$ 
The resulting matrix has size 
$$\bigl(\omega+(\infty+\infty)\bigr) \times \bigl(\omega+(\infty+\infty)\bigr) \quad=\quad
(\omega+\infty)\times(\omega+\infty)
,$$
i.e., we again get an element of $\Coll(\omega)$.

\begin{proposition}
 The product $\circ$ is a well-defined associative operation on the set $\Coll(\omega)$.
\end{proposition}

This can be verified by a straightforward calculation. For a clarification of this operation, 
see
\cite{Ner-book}, Section  IX.5.
 Classical operator colligations are matrices determined up to the equivalence
$$
\begin{pmatrix}\alpha&\beta\\ \gamma&\delta\end{pmatrix}
\sim
\begin{pmatrix}1&0\\0&u \end{pmatrix}
\begin{pmatrix}\alpha&\beta\\ \gamma&\delta\end{pmatrix}
\begin{pmatrix}1&0\\0&u^{-1} \end{pmatrix}
.
$$ 
Colligations, their multiplication, and characteristic functions 
appeared in the spectral 
theory
of non-self-adjoint operators (M.~S.~Livshits, V.~P.~Potapov, 1946--1955, \cite{Liv1},
 \cite{Liv2},
\cite{Pot}, see survey in \cite{Bro}, see also algebraic version in \cite{Dym}).


\smallskip

{\bf\punct Results of the paper.}
First (Theorem \ref{th:existence}),
 we prove the following statements:

\smallskip

--- The closure of $\GLO(\infty)$ in polymorphisms of $(\R^\infty,\mu_\infty)$  contains the semigroup 
$\Coll(\infty)$.

\smallskip

--- For $n<\infty$ the semigroup $\Coll(n)$ admits a canonical embedding to 
semigroup of polymorphisms of the space $(\R^n,\mu_n)$. 

\sm

Our main purpose 
is to write explicit formulas (Theorems  \ref{th:formula-1}, \ref{th:formula-2}) for this embedding. 

\smallskip


{\bf\punct A general problem.} Many interesting actions of infinite dimensional groups
on spaces with measures are known,  see survey \cite{Ner-frac} and recent 'new' constructions
\cite{KOV}, \cite{Pick}, \cite{Ner-hua}, \cite{BO}. In all  cases there arises the problem of 
description of closure of the group in polymorphisms,
in all the cases this gives  semigroups that  essentially differ from the initial groups%
\footnote{This is counterpart of Olshanski problem about weak closure of image of unitary
representation, see \cite{OlshGB}; for a finite-dimensional counterpart, see \cite{DCP}.} .
 In this work and in \cite{Ner-match} the problem was solved
in two the most simple cases (Gaussian and Poisson measures).
 In both cases we get unusual interesting formulas.  


\section{Polymorphisms. Preliminaries}

\COUNTERS

First, we need some preliminaries on polymorphisms.

{\bf\punct Measures on $\R^\times$.}
 Denote by $\R^\times$ the multiplicative group of positive real numbers,
denote by $t$ the coordinate on $\R^\times$, by $\phi*\psi$ we denote {\it convolution} of
finite measures $\phi$ and $\psi$ on $\R^\times$,  it  defined by
$$
\int_{\R^\times} f(t)\,d(\phi*\psi)(t)=
\int_{\R^\times}\int_{\R^\times} f(pq)\,d\psi(p)\,d\phi(q).
$$
Recall that a sequence of finite measures $\psi_j$ on $\R^\times$ 
{\it weakly converges} to a measure $\psi$ if for any 
continuous function $f$ on $\R^\times$ we have the convergence
$$
\int_{\R^"} f(t)\,d\psi_j(t)\longrightarrow\int_{\R^\times} f(t)\,d\psi(t)
.$$

\smallskip


{\bf\punct Product of polymorphisms.} Here we give a formal definition 
of the product of polymorphisms, but actially we use Theorem \ref{th:rep-cat}
instead of the definiton. For details, see \cite{Ner-poli}.

Let $p$ be a function on $M\times N$ taking values in finite measures on
$\R^\times$. Such a function determines a measure $\frP$ on a product
$M\times N\times \R^\times$,
$$
\iiint\limits_{M\times N\times \R^\times} f(m,n,t)\,d\frP(m,n,t)
:=\iint\limits_{A\times B}\int\limits_{\R^\times} f(m,n,t)\,dp(m,n)(t)\,d\nu(n)\,d\mu(m)
.$$
If $p$ satisfies two identities
\begin{align*}
\int_A \int_N \int_{\R^\times} d p(m,n)(t)\,dp(m,n)(t)\,d\nu(n)\,d\mu(m)=\mu(A),
\\
\int_M \int_B \int_{\R^\times}t\, d p(m,n)(t)\,\,dp(m,n)(t) \,d\nu(n)\,d\mu(m)=\nu(B)
\end{align*}
for any measurable subsets $A\subset M$, $B\subset N$,
then $\frP$ is a polymorphism. If $\frP$ has such aform,
we say that $\frP$ is absolutely continuous.

Now let $\frP\in\Pol(M,N)$, $\frQ\in\Pol(N,K)$ be absolutely continuous 
polymorphisms, $p$, $q$ be the correspondin functions. Then the function
$r$ on $M\times K$ is determined by
$$
r(a,c)=\int_N p(m,n)* q(n,k)\,d\nu(n).
$$
The integral is convergent a.s.

\begin{theorem}
 This product admits a unique separately continuous extension to an operation
$\Pol(M, N)\times\Pol(N,K)\to\Pol(M,K)$.
\end{theorem}

\sm


{\bf\punct Involution in the category of polymorphisms.}
Let $\frP:M\zigzag N$ be a polymorphism. We define the polymorphism
$\frP^\star:N\zigzag M$ by
$$
\frP^\star(n,m,t)=t\cdot\frP(m,n,t^{-1})
$$
For any polymorphisms $\frP:M\zigzag N$, $\frQ:N\zigzag K$,
the following property holds
$$
(\frQ\diamond\frP)^\star=\frP^\star\diamond\frQ^\star.
$$
If $g\in\Gms(M)$, then
$$
\frI(g)^\star=\frI(g^{-1})
.$$
 Our next purpose is to extend the operators (\ref{eq:T-lambda})  to arbitrary polymorphisms.


\smallskip

{\bf \punct Mellin transform  of polymorphisms.}
Here we present without proof some simple statements from \cite{Ner-poli}.
Notice that below we use Theorem \ref{th:rep-cat} and do not refer to the definition of product
of polymorphisms.

Fix $\lambda=\frac 1p+is\in \C$ as above (\ref{eq:lambda}). Let $q$
is defined from
$\frac 1p+\frac 1q=1$.
 For a polymorphism $\frP:M\zigzag N$ we consider the bilinear form
on $L^p(M,\mu)\times L^q(N,\nu)\to \C$ given by
$$
S_\lambda(f,g)=\iiint_{M\times N\times \R^\times} f(m) g(n) t^\lambda\,d\frP(m,n,t).
$$

\begin{proposition}
\label{pr:mellin-le}
 {\rm( \cite{Ner-poli})}
 {\rm a)}
$$
|S_\lambda(f,g)|\le \|f\|_{L_p}\cdot \|g\|_{L_q}.
$$

{\rm b)} $\frP$ is uniquely determined by the family of forms $S_\lambda(\cdot,\cdot)$.
\end{proposition}

\begin{corollary}
\label{cor:operator-exists}
 {\rm a} There exists a unique linear operator
$$
T_\lambda(\frP):L^p(N,\nu)\to L^p(M,\mu)
$$
such that
$$
S(f,g)=\int_M  f(m)\cdot T_\lambda(\frP)\cdot g(m)\,d\mu(m). 
$$

b)  $\|T_\lambda(\frP)\|\le 1$,  where a norm is the norm of an operator
$L^p(N,\nu)\to L^p(M,\mu)$.

\sm

{\rm c)} A polymorphism $\frP$ is uniquely determined
by the operator-valued function $\lambda\mapsto T_\lambda(\frP)$,
and, moreover, by its values on each line $\frac 1p+is$ for fixed $p$.
\end{corollary}

 For $h\in\Gms(M)$, we have
$$
T_\lambda(\iota(h))=T_\lambda(h)
,$$
where $T_\lambda(h)$ is defined by (\ref{eq:T-lambda}).

\begin{theorem}
\label{th:rep-cat}
 $T_\lambda$ is a representation of  a category, i.e.
\begin{equation}
T_\lambda(\frQ\diamond\frP)=
T_\lambda(\frQ)T_\lambda(\frP)
\label{eq:representation}
.\end{equation}
\end{theorem}


{\bf\punct Convergence.}

\begin{theorem}
\label{th:convergence}
{\rm a)} $T_\lambda(\frP)$ is weakly continuous, i.e., if $\frP_j$ converges to $\frP$, then
\begin{equation}
\int_M f(m)\cdot T_\lambda(\frP_j) g(m))\,d\mu(m) 
\quad \text{converges to}\quad
\int_M f(m) T_\lambda(\frP) g(m)\,d\mu(m) 
\label{eq:two-conver}
\end{equation}
for any $f\in L^q(M)$, $g\in L^p(N)$.

\smallskip

{\rm b)} Conversely, if {\rm (\ref{eq:two-conver})}
 holds for each $\lambda$ in the strip $0\le\Re \lambda\le 1$, then
$\frP_j$ converges to $\frP$. Moreover, it is sufficient
to require the convergences on the lines $\Re \lambda=0$ and $\Re\lambda=1$.
\end{theorem}

\section{Abstract statement}

\COUNTERS


\smallskip

{\bf \punct Polymorphisms $\frl_n$.} 
Let $(M,\mu)$ be a space with measure. Denote by $\Delta(m,m')$  the measure
on $M\times M$
supported by the  diagonal of $M\times M$ such that   the projection of $\Delta$ to
the first factor 
$M$ is $\mu$.

Let $\omega=0$, 1, \dots, $\infty$. 
Consider the space
$\R^\omega\times \R^\infty$ equipped with the measure $\mu_{\omega+\infty}=\mu_\omega\times \mu_\infty$.
Let $x$, $x'$ range in $\R^\omega$, $y$ in $\R^\infty$, $t$ in $\R^\times$.
Consider the polymorphism 
$$
\frl_\omega:(\R^\omega,\mu_\omega)\zigzag (\R^\omega\times \R^\infty,\mu_\omega\times\mu_\infty)
$$
given by
$$
\frl_\omega(x';x,y;t)=\Delta(x,x')\times \mu_\infty(y)\times \delta(t-1)
,$$
where $\delta$ is the delta-function.

The following  statement is straightforward.

\begin{lemma}
\label{l:Theta}
{\rm a)} For a function $f$ on $\R^\omega$ we have
 $$
T_\lambda(\frl_\omega)f(x,y)=f(x)
$$

{\rm b)} For a function $g(x,y)$ on $\R^{\omega+\infty}$, we have
 $$
T_\lambda(\frl_\omega^\star)g(x)=\int_{\R^\infty} g(x,y)\,d\mu_\infty(y)
$$ 

{\rm c)} $\frl_\omega^\star\diamond\frl_\omega:\R^\omega\zigzag\R^\omega$ is $\Delta(x,x')\times\delta(t-1)$.

\smallskip

{\rm d)} The polymorphism
 $$
\frt_\omega:= \frl_\omega\diamond\frl_\omega^\star:\R^{\omega+\infty}\zigzag\R^{\omega+\infty}
 $$
  equals
$$
\Delta(x,x')\times \mu_\infty(y)\times\mu_\infty(y')\times \delta(t-1)
,$$
where $(x,y)$ is  in the  first copy of $\R^{\omega+\infty}$ and $(x',y')$ is in the
second copy.

\smallskip

{\rm e)} The operator corresponding to $\frt_\omega$ is
$$
T_\lambda(\frt_\omega)f(x,y)=\int_{\R^\infty} f(x,z)\,d\mu_\infty(z).
$$
In particular, in $L^2$ this operator is the orthogonal projection to the space of functions 
independent on $y$.

\sm

{\rm f)} Consider a sequence $h_j=\begin{pmatrix}1&0\\0&u_j \end{pmatrix}\in\OO(\infty)$
where $u_j$ weakly converges to $0$. Then   $\frI(h_j)$ converges
to  $\frt_\omega=\frl_\omega\diamond\frl_\omega^\star$.
\end{lemma} 

An example of a sequence $u_j$ is 
$$
u_j=\begin{pmatrix} 0&1&0\\
                               1&0&0\\
                                0&0&1
\end{pmatrix}
\begin{array}{l}
\}j\\
\}j\\
\}\infty
\end{array}
$$


{\bf\punct Action of colligations.} Let $\omega=0$, 1, \dots, $\infty$.
Let $\fra\in\Coll(\omega)$, let $A$ be its representative in
$\GLO(\omega+\infty)$. Consider the polymorphism
$$
\tau^{(\omega)} (\fra):(\R^\omega,\mu_\omega)\zigzag (\R^\omega,\mu_\omega)
$$
given by
$$
\tau^{(\omega)}(\fra)=\frl_\omega \frI(A) \frl_\omega^\star.
$$

\begin{theorem}
\label{th:existence}
The map $\tau^{(\omega)}:\Coll(\omega)\to\Pol(\R^\omega,\R^\omega)$ is a homorphism of semigroups.
\end{theorem}

\begin{theorem}
\label{th:closure}
 For $\omega=\infty$ the image $\tau^{(\infty)}(\Coll(\infty))\subset \Pol(\R^\infty,\R^\infty)$ 
is contained in the closure
of $\frI\bigl(\GLO(\infty)\bigr)$.
\end{theorem}
%

{\bf\punct Proof of Theorem \ref{th:existence}.}
We must verify the identity
\begin{equation}
T_{\lambda}(\fra_1)T_{\lambda}(\fra_2)=T_{\lambda}(\fra_1\circ \fra_2)
.
\label{eq:hotim}
\end{equation}
or, equivalently,
$$
T_{\lambda}(\frt_\omega A_1\frt_\omega)T_{\lambda}(\frt_\omega A_2\frt_\omega)=
T_{\lambda}^{(\omega)}(\frt_\omega A_1A_2\frt_\omega)
.$$

Let $\rho$ be a unitary representation of $\GLO(\omega+\infty)\simeq\GLO(\infty)$ continuous with respect
 to the Shale topology.
Denote by $H(\omega)$ the space of $\OO(\infty)$-invariant vectors. 
Denote by $P(\omega)$ the orthogonal projection
on $H(\omega)$.  For $A\in\GLO(\omega+\infty)$, we define the operator
\begin{equation}
\rho^{(\omega)}(\fra):=P(\omega)\rho(A):\,\, H(\omega)\to H(\omega).
\label{eq:srezka}
\end{equation}
It can be easily checked that $\rho^{(\fra)}(g)$ depends on a operator colligation $\fra$
 and not on $A$ itself. 

\begin{theorem} 
 We get a representation of the semigroup $\Coll(\omega)$ in the space $H(\omega)$.
\begin{equation}
\rho^{(\omega)}(\fra_1)\rho^{(\omega)}(\fra_2)=\rho^{(\omega)}(\fra_1\circ \fra_2)
.
\label{eq:multiplicativity}
\end{equation}

\end{theorem}

See \cite{OlshGB}, \cite{Ner-book}, see a simple proof in \cite{Ner-spheric}.

We need this theorem for representations $T_{1/2+is}$ of the group $\GLO(\omega+\infty)$
in $L^2(\R^{\omega+\infty}),\mu_{\omega+\infty}$, in this case $P(\omega)$
is $T_{1/2+is}(\frt)$, 
$$
T_{1/2+is}(\fra)=T_{1/2+is}(\frt)T_{1/2+is}(A)T_{1/2+is}(\frt)
,$$
the identity \ref{eq:multiplicativity} can be written as
\begin{equation}
T_{1/2+is}^{(\omega)}(\fra_1)T_{1/2+is}^{(\omega)}(\fra_2)=T_{1/2+is}^{(\omega)}(\fra_1\circ \fra_2)
\label{eq:olsh}
\end{equation}

Since $T_\lambda$ depends holomorphically in $\lambda$, we get
(\ref{eq:hotim}).

\sm

{\sc Remark.} Identity \ref{eq:olsh} can be verified by a long straightforward calculation
(and in fact this was done in \cite{OlshGB}).

\smallskip



{\bf\punct Proof of Theorem \ref{th:closure}.} 
Let $\fra\in\Coll(\infty)$, let $A\in\GLO(\infty+\infty)$
be its representative. We define the polymorphism
$$
\sigma(\fra):(\R^{\infty+\infty},\mu_{\infty+\infty})\zigzag (\R^{\infty+\infty},\mu_{\infty+\infty})
$$
by 
$$
\sigma (\fra)=\frt_\infty \diamond \tau(A)\diamond
 \frt_\infty^\star
 .
$$
By Lemma \ref{l:Theta}.f, the element $\frt_\infty$ is contained in the closure
of $\O(\infty)$.  By separate continuity of the product, 
$\frt_\infty \diamond \tau(A)\diamond \frt_\infty^\star$ is contained in the closure
of $\GLO(\infty+\infty)$

Next, represent the set of natural numbers $\N$ as a union of two disjoint sets $I$, $J$.
Consider the monotonic  bijections $I\to\N$,  $J\to\N$. In this way
 we identify $\R^\infty$ and $\R^{\infty+\infty}$.
Denote by $\sigma(\fra; I):\R^\infty \zigzag \R^\infty$
 the image of the polymorphism $\sigma(\fra)$ under this identification.
By construction $\sigma(\fra,I)$ is contained in the closure of $\GLO(\infty)$.

Now take 
$$
I_k=\{\text{1, 2, 3, \dots, $k$, $k+2$, $k+4$, $k+6$,\dots}\},
$$
 Then $\sigma(\fra,I_k)$ converges to $\tau(\fra)$.
\hfill $\square$


\smallskip

{\bf \punct Injectivity.}
We formulate without proof the following statement.

\begin{theorem}
The maps $\Coll(\omega)\to \Pol(\R^\omega,\R^\omega)$ are injective.
\end{theorem}

 This is equivalent to the statement: the family of representations 
 $\fra\mapsto P(\omega)T_\lambda(\fra)P(\omega)$ 
 separates points of $\Coll(\omega)$. 


\section{Canonical forms}

\COUNTERS

{\bf\punct  Canonical forms.} Let $n<\infty$, $\frg\in\Coll(n)$.  
Let $g=\begin{pmatrix} g_{11}&g_{12}\\g_{21}&g_{22}\end{pmatrix}$ be a representative of $\frg$. 

\begin{lemma}
\label{l:canonical}
Assume that rank of $g_{12}$ is maximal. Then $\frg$ has a representative of the form
\begin{equation}
\begin{array}{c}G=
\\ \phantom{=}
\\ \phantom{=}
\end{array}
\begin{array}{cc}
\begin{pmatrix}
a&b\\c&d\\0&H
\end{pmatrix}
\begin{array}{l}
\}n\\\}n\\\}\infty
\end{array}
\\
\begin{matrix}\!\!\!
\small \underbrace{\hphantom{1}}_{n}&\!\!\!\!\!\!
\small \underbrace{\hphantom{1}}_{n+\infty}
& \hphantom{\}\infty}
\end{matrix}
\end{array}
\quad
\begin{matrix}
=\\\phantom{=}
\\\phantom{=}
\end{matrix}
\quad
\begin{array}{cc}
\begin{pmatrix}a&b_1&b_2\\
c&d_1&d_2\\
0&0&h
\end{pmatrix}
\begin{array}{l}
\}n\\\}n\\\}\infty
\end{array}
\\
\begin{matrix}\!\!\!\!\!\!\small\underbrace{\hphantom{1}}_n
&\!\!\!\!\!\!
\small \underbrace{\hphantom{1}}_{n}&\!\!\!\!\!\!
\small \underbrace{\hphantom{1}}_{\infty}
& \hphantom{\}\infty}
\end{matrix}
\end{array}
\label{eq:can-form}
\end{equation}
where $h$ is a diagonal matrxix with positive entries $h_j$, $\sum(h_j-1)^2<\infty$.
\end{lemma}

\begin{lemma}
\label{l:oo}
Any $g=\begin{pmatrix}\alpha &\beta\\ \gamma&\delta\end{pmatrix} \in \O(n+\infty)$
admits a representation in the form 
$$
g=(1+S)\begin{pmatrix}1& 0\\ 0& u \end{pmatrix}
,$$
where $S$ is a Hilbert--Schmidt matrix and $u\in \O(\infty)$.
\end{lemma}

{\sc Proof of Lemma \ref{l:oo}.} The matrix $\delta^t\delta-1$ 
is Hilbert--Schmidt and $\delta$ is Fredholm of index 0, therefore $\delta$ can be represented as
$$
\delta= v H u
,$$
where $u$, $v\in\O(\infty)$, and $H$ is a diagonal matrix, the matrix $H-1$ is Hilbert--Schmidt.
Therefore $g$ has the form 
$$
g= \begin{pmatrix} 1&0\\0&v \end{pmatrix}
 \begin{pmatrix}\alpha&\beta' \\ \gamma'& H\end{pmatrix}
\begin{pmatrix} 1&0\\0&u \end{pmatrix}
$$
The middle factor is ($1+$ Hilbert--Schmidt matrix).  Finally, we get a desired representation
$$
g= \left[\begin{pmatrix} 1&0\\0&v \end{pmatrix}
 \begin{pmatrix}\alpha&\beta' \\ \gamma'& H\end{pmatrix}
 \begin{pmatrix} 1&0\\0&v \end{pmatrix}^{-1}\right]
 \cdot
\left[ \begin{pmatrix} 1&0\\0&v \end{pmatrix}
\begin{pmatrix} 1&0\\0&u \end{pmatrix}
\right]
$$

{\sc Proof of Lemma \ref{l:canonical}.}
By Lemma \ref{l:oo}, we can assume that $G-1$ is a Hilbert--Schmidt
matrix.
Since $\rk g_{12}=n$, a left multiplication by an orthogonal matrix $w$ can reduce $g_{12}$
to the form $\begin{pmatrix}c\\0\end{pmatrix}$.

 Thus we get a matrix
$R'=\begin{pmatrix} a&b\\c&d\\0&H\end{pmatrix}$ such that $R'-1$ is Hilbert--Schmidt.
 We transform $R'$ by 
$$
\begin{pmatrix} a&b\\c&d\\0&H\end{pmatrix}
\longrightarrow
\begin{pmatrix}
1&0&0\\
0&1&0\\
0&0&u
\end{pmatrix}
\begin{pmatrix} a&b\\c&d\\0&H\end{pmatrix}
\begin{pmatrix}
1&0&0\\
0&v_{11}&v_{12}\\
0&v_{21}&v_{22}
\end{pmatrix},
$$
where $u$ and $\begin{pmatrix}
v_{11}&v_{12}\\v_{21}&v_{22}
\end{pmatrix}
$ are orthogonal matrices. Consider 
$(n+\infty)\times\infty$ matrix $J=\begin{pmatrix}0 & 1\end{pmatrix}$.
Then $H-J$ is a Hilbert--Shmidt operator,
therefore the Fredholm index of $H$ equals $n$. 
 Since $G$ is invertible,  $\ker H=0$,
Hence
 $\mathrm{codim}\, \Im H=n$. Such $H$ can be reduced to the form
$\begin{pmatrix}0&h \end{pmatrix}$, where $h$ is diagonal.
The standard proof of the theorem about singular values (see \cite{RS1})
 can be adapted to this case. 
\hfill $\square$

\smallskip

{\bf\punct  Coordinates.}
Take a colligation reduced to a canonical form (\ref{eq:can-form}). We pass to
{\it Potapov coordinates} (see \cite{Pot}) on the space of matrices,
$$
\begin{pmatrix} P&Q\\R&T\end{pmatrix}
:=\begin{pmatrix}b-ac^{-1}d&- ac^{-1}\\ c^{-1}d&c^{-1}\end{pmatrix}
$$
or
$$
\begin{pmatrix} P_1&P_2&Q\\R_1&R_2&T\end{pmatrix}
:=\begin{pmatrix}b_1-ac^{-1}d_1&b_2-ac^{-1}d_2&- ac^{-1}\\ c^{-1}d_1&c^{-1}d_2&c^{-1}\end{pmatrix}
,$$
the size of the block matrices is $(n+\infty+n)\times(n+n)$.
Formulas below are written in the terms of $P$, $Q$, $R$, $T$, and $h$.

\section{Calculations. Finite matrices}

\COUNTERS

{\bf\punct Measures  $\Phi[b,M;t]$.} Let $M\ge 0$, $b\in\R$.
We define the measure $\Phi[b,M;t]$ on $\R^\times$ by

\sm

--- for $b>0$
$$
\Phi[b,M;t]=
\begin{cases}
\frac 1{\sqrt{2\pi}}
t^{1/b}
(- b\ln t)^{-1/2}\cosh\sqrt{-\frac{4M}b\ln t} \,\frac{dt}t & \text{if $0<t<1$};
\\
0 & \text{if $t>1$}.
\end{cases}
$$
 
--- for $b=0$
$$
\Phi[0,M;t]=e^{M}\delta(t-1)
$$
 
--- for $b<0$,
$$
\Phi[b,M;t]=
\begin{cases}
0&\text{if $0<t<1$}
\\
\frac 1{\sqrt{2\pi}}
t^{-1/b}(4M b\ln t)^{-1/2}\cosh\sqrt{\frac{4M}b\ln t}\,\frac{dt}t &
\text{if $t>1$}
\end{cases}
$$

\begin{lemma}
\label{l:Phi}
$$
\frac 1{\sqrt{2\pi}}\int_{\R^\times} t^\lambda \Phi[b,M;t]=\frac1{\sqrt{1+b\lambda}}
\exp\Bigl\{ \frac{M}{1+b\lambda} \Bigr\}
.
$$
\end{lemma}

{\sc Proof.} To be definite, set $b>0$. We must evaluate
$$
\frac 1{\sqrt{2\pi}}
\int_0^1 t^{\lambda+1/b}(- b\ln t)^{-1/2}
\cosh
\sqrt{-\frac{4M}b\ln t}
\,\frac {dt}t
.
$$
We substitute $y=\ln t $ and get
$$
\frac 1{\sqrt{2\pi}}
\int_{-\infty}^0
e^{(\lambda+1/b)y}(- by)^{-1/2}
\cosh\sqrt{-\frac{4M}b y}\,dy
.
$$
Next, we set $z=-\frac{4M}b y$, and come to
\begin{multline*}
\frac 1{\sqrt{2\pi}\cdot\sqrt{4M}}
\int_0^\infty e^{-\frac 1{4M}(b\lambda+1)z}z^{-1/2}\cosh \sqrt z\, dz 
=\\=
\frac 1{\sqrt{2\pi}\cdot\sqrt{M}}
\int_0^\infty e^{-\frac 1{4M}(b\lambda+1)u^2}\cosh u\,du
.
\end{multline*}
Writing $\cosh u=\frac 12(e^u+e^{-u})$, we get
$$
\frac 1{\sqrt{2\pi}\cdot 2\sqrt{M}}
\int_{-\infty}^\infty e^{-\frac 1{4M}(b\lambda+1)u^2}e^u\,d u=
\frac 1{\sqrt{1+b\lambda}}
\exp\Bigl\{\frac {M}{1+b\lambda}\Bigr\}
.
$$


\smallskip


{\bf \punct Formula.} We consider coordinates on $\Coll(n)$ defined above.
For $x$, $u\in\R^n$ we define the following $\delta$-measure
$dN_{x,u}(t)$ on $\R^\times$
$$dN_{x,u}(t)=A(x,u)\,\delta\bigl(t-B(x,u)\bigr),$$
where
$$
A(x,u)=
|\det T|
\exp\Bigl\{ -\frac12\|xQ+uT\|^2-\frac 12 \|(xP+uR) H^t(1- HH^t)^{-1}\|^2    \Bigr\} 
,
$$
\begin{multline}
B(x,u)=|\det G|
\exp\Bigl\{
 \frac 12\bigl(\|xQ+uT\|^2 -\|x\|^2+\|u\|^2-\\- (xP+uR)(1-H^tH)^{-1} (xP+uR)^t\bigr)\Bigr\}
,\end{multline}
where $\|\cdot\|$ is the standard norm in $\R^n$.

Denote by $h_j$ the diagonal entries of the matrix $h$. Denote by
$(\psi_1,\psi_2,\dots)$ the coordinates of the vector $xP_2+uR_2$.

\begin{theorem}
\label{th:formula-1}
Let $\frg\in\Coll(n)$ have a representative 
\begin{equation}
G=
\begin{array}{cc}
\begin{pmatrix}a&b_1&b_2&0\\
c&d_1&d_2&0\\
0&0&h&0\\
0&0&0&1\\
\end{pmatrix}
\begin{array}{l}
\}n\\\}n\\\}m-n\\\} \infty
\end{array}
\\
\begin{matrix}\!\!\!\!\!\!\small\underbrace{\hphantom{1}}_n
&\!\!\!\!\!\!
\small \underbrace{\hphantom{1}}_{n}&
\!\!\!\!\!\!\small \underbrace{\hphantom{1}}_{m-n}&
\!\!\!\!\!\!\small \underbrace{\hphantom{1}}_{\infty}
& \hphantom{\}\infty}
\end{matrix}
\end{array}
\label{eq:G-finite}
\end{equation}
and $h_j\ne 1$.
  Then the polymorphism $\tau(\fra)$ is given by
\begin{equation}
\biggl(N_{x,u}(t)* 
\begin{array}{c}m-n\\\text{\Huge $\ast$}\\ j=1\end{array}
 \quad \Phi\Bigl[h_j^2-1, \frac{h_j^2|\psi_j|^2}{2(1-h_j^2)};\,t\Bigr]\biggr)dx\,du
,\
\label{eq:final-1}
\end{equation}
where $*$ denotes the convolution in $\R^\times$ and {\Huge $\ast$} is the symbol of multiple convolution
with respect to $j$.
\end{theorem}


{\bf\punct Transformation of the determinant.}
 Note that
\begin{multline*}
\det G=
\det\begin{pmatrix}
a&b_1&b_2\\
c&d_1&d_2\\
0&0&h
\end{pmatrix}
=\\=
\det\begin{pmatrix}
a&b_1\\
c&d_1
\end{pmatrix}\cdot \det(h)=
\pm
\det(c)\det(b_1-a c^{-1} d_1)\det(h)
.
\end{multline*}
Thus
$$
|\det G|
=
\left|\frac{\det (P_1)\det(H)}{\det (T)}\right|
.$$

{\bf\punct Calculation.}
 We wish to write explicitly operators
(\ref{eq:srezka}) for the representations $T_\lambda(G)$. 
$$
T_\lambda^{(n)}(G)=
T_\lambda(\frl) T_\lambda(G) T_\lambda(\frl^\star)
.
$$
Let $x\in\R^n$, $y\in \R^n$, $z\in\R^{m-n}$, $\xi\in\R^\infty$.
The operator $T_\lambda(\frl^\star)$ sends  
a function $f(x)$ on $\R^n$ to the same function
$f(x)$ on $\R^n\times \R^n \times \R^{m-n}\times \R^\infty$.
We apply $T_\lambda(G)$ and come to
\begin{equation}
|\det G|^\lambda
 f(xa+yc)
\exp\left\{-\frac \lambda 2 \begin{pmatrix} x&y&z\end{pmatrix}
(GG^t -1)\begin{pmatrix} x^t\\y^t\\z^t\end{pmatrix}\right\}
.
\label{eq:prom1}
\end{equation}
Next, the operator $T_\lambda(\frl)$ is  the average with respect to 
variables 
$(y,z,\xi)\in \R^n \times \R^{m-n}\times \R^\infty$.
Since the function (\ref{eq:prom1})
is independent on $\xi$, we take average with respect to
$(y,z)$. We come to
 \begin{multline}
T_\lambda^{(n)}(G)f(x)=
|\det G|^\lambda
\iint\limits_{\R^n\times \R^{m-n}}
 f(xa+yc)\times\\ \times
\exp\left\{-\frac \lambda 2 \begin{pmatrix} x&y&z\end{pmatrix}
(GG^t -1)\begin{pmatrix} x^t\\y^t\\z^t\end{pmatrix}\right\}
\,d\mu_n(y)\, d\mu_{m-n}(z)
=\\=
\frac{|\det(G)|^\lambda}
{(2\pi)^{m/2}}
\cdot
e^{\frac1 2 x^2}
\iint\limits_{\R^n\times \R^{m-n}} 
f(xa+yc)\times\\ \times
\exp\left\{-\frac \lambda 2 \begin{pmatrix} x&y&z\end{pmatrix}
GG^t \begin{pmatrix} x^t\\y^t\\z^t\end{pmatrix}
+
\frac {\lambda-1} 2 \begin{pmatrix} x&y&z\end{pmatrix}
\begin{pmatrix} x^t\\y^t\\z^t\end{pmatrix}
\right\}
\,dy\,dz
\label{eq:int-operator}
\end{multline}

We change variable $y$ by $u$ according
$$
u=xa+yc,\qquad y=uc^{-1}-xac^{-1}
.
$$
Then
$$
\begin{pmatrix}x&y&z\end{pmatrix}=\begin{pmatrix} x&u&z\end{pmatrix} S
,$$
where 
$$
S=\begin{pmatrix}1&-ac^{-1}&0\\
                         0&c^{-1}&0\\
                         0&0&1
 \end{pmatrix}
.
$$
Quadratic form in (\ref{eq:int-operator}) transforms to
$$
\left\{-\frac \lambda 2 \begin{pmatrix} x&u&z\end{pmatrix}
SGG^tS^t \begin{pmatrix} x^t\\u^t\\z^t\end{pmatrix}
+
\frac {\lambda-1} 2 \begin{pmatrix} x&u&z\end{pmatrix}
SS^t
\begin{pmatrix} x^t\\u^t\\z^t\end{pmatrix}
\right\}
$$
Passing to Potapov coordinates, we get
$$
SS^t=\begin{pmatrix}1+QQ^t&QT^t&0\\ TQ^t&TT^t&0\\0&0&1 \end{pmatrix}
$$
$$
SG=\begin{pmatrix}0&P\\1&R\\0&H \end{pmatrix}
\qquad
SGG^tS^t=\begin{pmatrix} PP^t& PR^t& PH^t\\
RP^t& 1+RR^t&RH^t\\
HP^t&HR^t&HH^t
\end{pmatrix}
$$
We come to the expression of the form
$$
T_\lambda^{(n)}(G)\,f(x)=
\int_{\R^n} \cK(x,u) f(u)\, du
,$$
where
 the kernel $\cK$ is given by
$$
\cK(x,u)=(2\pi)^{-n/2}
|\det(G)|^\lambda |\det c|^{-1} \exp\{V(x,u)\}
 \int_{\R^{m-n}}\exp\bigl\{ U(x,u,z)\bigr\} \,dz
,$$
where
\begin{multline}
\exp\bigl\{V(x,u)\bigr\}=\exp\Bigl\{\frac 1 2 x x^t+\frac{\lambda-1}2
\begin{pmatrix}x&u \end{pmatrix}
\begin{pmatrix}
QQ^t+1&QT^t\\TQ^t&TT^t
\end{pmatrix}
\begin{pmatrix}x^t\\u^t \end{pmatrix}
-\\-
\frac\lambda2
\begin{pmatrix}x&u \end{pmatrix}
\begin{pmatrix}PP^t&PR^t\\RP^t&RR^t+1\end{pmatrix}
\begin{pmatrix}x^t\\u^t\end{pmatrix}
\Bigr\}
=\\=
\exp\Bigl\{-\frac\lambda 2\|xP+uR\|^2+\frac{\lambda-1}2 \|xQ+uT\|^2+
\frac\lambda 2(\|x\|^2-\|u\|^2)\Bigr\}
\label{eq:1}
\end{multline}
and 
\begin{multline}
\int_{\R^{m-n}}\exp\Bigl\{U(x,u,z)\Bigr\}\,dz=\\=
{(2\pi)}^{-(m-n)/2}
\int_{\R^{m-n}}\exp\Bigl\{\frac 12 z(-\lambda HH^t+\lambda-1)z^t\Bigr\}
\exp\Bigl\{-\lambda zH (P^tx^t+R^tu^t)\Bigr\}\,dz=
\\=
\det(\lambda HH^t-\lambda+1)^{-1/2}
\times\\ \times 
\exp\Bigl\{\frac{\lambda^2} 2(xP+uR)H^t(\lambda HH^t-\lambda+1)^{-1}H(xP+yR)^t\Bigr\}
\label{eq:int-result}
\end{multline}

We wish to examine the exponential factor in (\ref{eq:int-result}).
Recall that $H$ is an $(m\times n)$ matrix of the form
$$
H=\begin{pmatrix} 0&\dots&0&h_1&0&\dots&0\\
                             0&\dots&0&0& h_2&\dots&0\\
                             \vdots&\ddots&\vdots&\vdots&\vdots&\ddots&\vdots\\
                             0&\dots&0&0&0&\dots&h_{m-n}
 \end{pmatrix}
$$
Therefore $HH^t$ is the diagonal matrix with entries
$h_j^2$ and $H^t(\lambda HH^t-\lambda+1)^{-1}H$ is the diagonal matrix with entries
 $0$ ($n$ times) and 
$
\frac{h_j^2}{\lambda h_j^2-\lambda+1}
.$
Therefore,
(\ref{eq:int-result}) equals 
\begin{equation}
(2\pi)^{n-m}
\prod_{j=1}^{m-n}
\bigl(1+\lambda ({h_j^2-1})\bigr)^{-1/2}
\exp\Bigl\{  \frac{\lambda^2 h_j^2 |\psi_j|^2}{2(\lambda h_j^2-\lambda+1)}\Bigr\}
\label{eq:product-1}
\end{equation}
Next, we write
\begin{equation}
\label{eq:divergence}
 \frac{\lambda^2 h_j^2}{\lambda h_j^2-\lambda+1}
=
\frac{\lambda h_j^2}{h_j^2-1}-\frac{h_j^2}{(h_j^2-1)^2}+ \frac{h_j^2}{(h_j^2-1)^2}\cdot
\frac 1{\lambda h_j^2-\lambda+1}
\end{equation}
and represent the product (\ref{eq:product-1}) as
\begin{multline}
\exp\Bigl\{-\frac 12 (xP+uR) H^t(1- HH^t)^{-2} H (xP+uR)^t \Bigr\}
\times\\ \times
\exp\Bigl\{-\frac \lambda 2(xP+uR)H^t(1-HH^t)^{-1} H (xP+uR)^t\Bigr\}
\times\\ \times
\prod_{j=1}^{m-n}
(\lambda (h_j^2-1)+1)^{-1/2}
\exp\Bigl\{\frac{h_j^2\|\psi_j\|^2}{2(h_j^2-1)^2}\cdot
\frac 1{\lambda(h_j^2-1)+1} \Bigr\}
\label{eq:2}
\end{multline}

Uniting (\ref{eq:1}) and (\ref{eq:2}), we come to a final expression for the kernel
of integral operator
\begin{align}
&\cK_\lambda(x,u)=\nonumber\\=
 & |\det c|^{-1}
\exp\Bigl\{ -\frac12\|xQ+uT\|^2-\frac 12 \|(xP+uR) H^t(1- HH^t)^{-1}\|^2    \Bigr\} 
\times
\label{eq:line-1} \\ 
&\times
|\det(G)|^\lambda \cdot
\exp\Bigl\{
\frac{\lambda}2\bigl( \|xQ+uT\|^2 +\|x\|^2-\|u\|^2-
\label{eq:line-2}
\\&\qquad\qquad\qquad\qquad- (xP+uR)(1-H^tH)^{-1} (xP+yR)^t
\bigr)\Bigr\}
\times
\label{eq:line-3}
\\ 
&\times
\prod_{j=1}^{m-n}
(\lambda (h_j^2-1)+1)^{-1/2}
\exp\Bigl\{\frac{h_j^2\|\psi_j\|^2}{2(h_j^2-1)^2}\cdot
\frac 1{\lambda(h_j^2-1)+1} \Bigr\}
.
\label{eq:line-4}
\end{align}

Now we must represent the kernel as a Mellin transform of a measure
$$
\cK_\lambda(x,u)=\int_0^\infty t^\lambda d M_{x,u}(t)
.
$$
The expression for $\cK_\lambda(x,u)$ is a product, therefore its Mellin
transform is a convolution. We must evaluate inverse Mellin transform for all factors.
The first factor (\ref{eq:line-1}) is  constant. The second factor 
(\ref{eq:line-2})--(\ref{eq:line-3}) has the form
$e^{\lambda a(x,u)}$, we have
$$
e^{\lambda a(x,u)}=\int_0^\infty t^\lambda \delta\bigl(t- e^{a(x,u)}\bigr)
.$$
For factors in (\ref{eq:line-4}) the inverse Mellin transform was evaluated
in Lemma \ref{l:Phi}.

This proves Theorem \ref{th:formula-1}.

\smallskip


\section{Convergent formula}

\COUNTERS

{\bf\punct Formula.} Now consider  arbitrary $\frg\in\Coll(n)$ being in the canonical form
(\ref{eq:can-form}),
$$
\begin{pmatrix}
a&b_1&b_2\\
c&d_1&d_2\\
0&0&h\\
\end{pmatrix}
$$
To write a formula that is valid in general case, we rearrange factors in (\ref{eq:final-1}).
First, we define $\delta$-measures on $\R^n\times\R^n$ by
$$dN^\circ_{x,u}(t)=A^\circ(x,u)\delta\bigl(t-B^\circ(x,u)\bigr),$$
where
$$
A^\circ(x,u)=
\det (T)\,
\exp\Bigl\{ -\frac 12\|xQ+uT\|^2    \Bigr\} 
$$
$$
B^\circ(x,u)=\frac{|\det P_1|}{|\det T|}
\exp\Bigl\{
 \frac 12\bigl(\|xQ+uT\|^2-\|xP_1+uR_1\|^2 -\|x\|^2+\|u\|^2\bigr)\Bigr\}
.$$
In fact, $dN^\circ_{x,u}(t)$ is the measure $dN_{x,u}(t)$ defined for the matrix
$\begin{pmatrix}a&b_1\\c&d_1\end{pmatrix}$.

Next, we define the following probability measures $\Xi_j=\Xi[h_j, \psi_j]$ 
on $\R^\times$:
\begin{multline}
\Xi[h_j, \psi_j]=\\=
\exp\Bigl\{- \frac{|\psi_j|^2 h_j^2}{2(1-h_j^2)^2}\Bigr\}
\cdot \delta\Bigl(t-h_j \exp\Bigl\{\frac{|\psi_j|^2}{2(1-h_j^2)}\Bigr\}\Bigr)
*\Phi\Bigl[h_j^2-1, \frac{h_j^2 |\psi_j|^2}{2(1-h_j^2)^2};\,t\Bigr]
\label{eq:Xi}
\end{multline}
if $h_j\ne 1$. For $h_j=1$ we set 
 $$
\Xi[1,\psi_j]=
\frac{1}{|\psi_j|}e^{-\frac 18|\psi_j|^2}
\exp\Bigl\{- \frac{\ln^2 t}{2|\psi_j|^2}\Bigr\}\frac {dt} {t^{3/2}},
\qquad \Xi[1,0]=\delta(t-1).$$

\begin{theorem}
\label{th:formula-2}
Let $\fra\in\Coll(n)$ be arbitrary. Then the polymorphism
$\tau(\fra)$ is given by
\begin{equation}
\Bigl(
dN^\circ_{x,u}(t)*
\begin{array}{c}\infty\\\text{\Huge $\ast$}\\ j=1\end{array} 
\Xi[h_j,\psi_j]\Bigr) dx\,du
.
\label{eq:final-2}
\end{equation}
\end{theorem}

\begin{lemma}
\label{l:fin1}
{\rm a)} Measures $\Xi[h_j,\psi_j]$ are probabilistic.

\sm

{\rm b)} The products
\begin{equation}
\begin{array}{c}\infty\\\text{\Huge $\ast$}\\ j=1\end{array} 
\Xi[h_j,\psi_j],
\qquad  \qquad
\begin{array}{c}\infty\\\text{\Huge $\ast$}\\ j=1\end{array} 
 \bigl( t\cdot \Xi[h_j,\psi_j]\bigr)
\label{eq:long-convolution}
\end{equation}
weakly converge in the semigroup of measures on $\R^\times$.
\end{lemma}

\begin{theorem}
\label{l:fin2}
{\rm a)}  For a matrix $g$ denote by 
denote by $g^{(m)}$
the matrix $\begin{pmatrix} z&0\\ 0&1\end{pmatrix}$, where $z$ is the upper left
$(n+m)\times(n+m)$ corner of the matrix $g$. Then the the polynorphism
$\tau(\frg^{(m)})$ coincides with
\begin{equation}
\Bigl(
dN^\circ_{x,u}(t)*
\begin{array}{c}m-n\\\text{\Huge $\ast$}\\ j=1\end{array} 
\Xi[h_j,\psi_j]\Bigr) dx\,du
\label{eq:final-1.5}
.
\end{equation}

{\rm b)} The sequence of polymorphisns
{\rm(\ref{eq:final-1.5})}
converges in semigroup of polymorphisms of $(\R^n,\mu_n)$.
 to $\tau(\fra)$.
\end{theorem}


{\bf\punct Rearrangement of factors (Lemma \ref{l:fin2}.a.}
First, rearrange factors in (\ref{eq:line-1})--(\ref{eq:line-4}):
\begin{align}
&\cK_\lambda(x,u)=
  |\det T|
\exp\Bigl\{ -\frac12\|xQ+uT\|^2\Bigr\} 
\Bigl(\frac{|\det(P_1)|}{\det (T)|}\Bigr)^\lambda
\times
\label{eq:line-5}
\\& \times
\exp\Bigl\{
\frac{\lambda}2\bigl( \|xQ+uT\|^2 +\|x\|^2-\|u\|^2-
\|xP_1+uR_1\|^2)\Bigr\}
\label{eq:line-6}
\\ 
&\times
\prod_{j=1}^{m-n}
\biggl(\exp\Bigl\{\frac{h_j^2 |\psi_j|^2}{2(1-h_j^2)^2} \Bigr\}
\cdot h_j^\lambda
\exp\Bigl\{\frac{\lambda|\psi_j|^2}{2(1-h_j^2)}
\Bigr\}
\times 
\label{eq:line-7}
\\
&\qquad\qquad\times
\bigl(\lambda (h_j^2-1)+1\bigr)^{-1/2}
\exp\Bigl\{\frac{h_j^2\|\psi_j\|^2}{2(h_j^2-1)^2}\cdot
\frac 1{\lambda(h_j^2-1)+1} \Bigr\}
\biggr)
\label{eq:line-8}
\end{align}

Factors in the product (\ref{eq:line-5})--(\ref{eq:line-6}) looks as singular near $h_j=1$.
But this singularity is artificial, it appears due division in the line (\ref{eq:divergence}).
Returning to the previous line (\ref{eq:product-1}) of the calculation, we get 
for $h_j=1$ the following factor
$$
\exp\Bigl\{-\frac12\lambda |\psi_j|^2+\frac12\lambda^2|\psi_j|^2 \Bigr\}
=\frac{1}{|\psi_j|}e^{-\frac 18|\psi_j|^2}\int_0^\infty t^\lambda \, 
\exp\Bigl\{- \frac{\ln^2 t}{2|\psi_j|^2}\Bigr\}\frac {dt} {t^{3/2}}
$$


{\bf \punct Proof of Lemma \ref{l:fin2}.b).} 

\begin{lemma}
The embedding $\iota:\GLO(\infty)\to \Pol(\R^\infty, \R^\infty)$ is continuous.
\end{lemma}

{\sc Proof.}  According Proposition \ref{th:convergence}.b it is sufficient to prove that 
the representations $T_\lambda(g)$ of $\GLO(\infty)$ are weakly continuous for all $\lambda$.
It is sufficient to take $f=e^{iax}$ and $g=e^{ibx}$ in  (\ref{eq:two-conver}) 
and to verify continuity of the corresponding matrix elements 
with respect to the Shale topology. \hfill$\square$

 \smallskip

Let $g$ be of the form (\ref{eq:can-form}).  
For finite matrices formulas (\ref{eq:final-1}) and (\ref{eq:final-2})
coincide.
Denote by $g^{(m)}$
the matrix $\begin{pmatrix} z&0\\ 0&1\end{pmatrix}$, where $z$ is the upper left
$(n+m)\times(n+m)$ corner of the matrix $g$.
For $g^{(m)}$ the formula (\ref{eq:final-1.5}) gives a correct result.
Next, $g^{(m)}$ converges to $g$ in the Shale topology. 
Therefore $\tau(g^{(m)})$ converges to $\tau(g)$ as $g\to\infty$.
This proves the last statement of the theorem.


\sm

%

\sm


{\bf \punct Proof of Theorem \ref{th:formula-2}.}
We must  prove convergence of the infinite convolution
in (\ref{eq:long-convolution}).
 The characteristic function of $\Xi[h_j,\psi_j]$ is given by
\begin{equation*}
\int_0^\infty t^\lambda \Xi_j[h_j,\psi_j]=
 h_j^\lambda
\bigl(1+\lambda ({h_j^2-1})\bigr)^{-1/2}
\exp\Bigl\{  \frac{\lambda^2 h_j^2 |\psi_j|^2}{2(\lambda h_j^2-\lambda+1)}-
\frac\lambda 2 |\psi_j|^2\Bigr\}
\end{equation*}
We have $\sum (h_j-1)^2<\infty$, $\sum|\psi_j|^2<\infty$. Under these conditions
we have a  convergence of the product in the strip $0\le \Re\lambda\le1$. This implies
the weak convergence of measures on $\R^\times$.

The convergence is uniform on compacts sets with respect to $x$, $u$, and this implies coincidence
of (\ref{eq:final-2}) and limit of (\ref{eq:final-1.5}).



{\tt Math.Dept., University of Vienna,

 Nordbergstrasse, 15,
Vienna, Austria

\&

Institute for Theoretical and Experimental Physics,

Bolshaya Cheremushkinskaya, 25, Moscow 117259,
Russia

\&

Mech.Math. Dept., Moscow State University,
Vorob'evy Gory, Moscow

e-mail: neretin(at) mccme.ru

URL:www.mat.univie.ac.at/$\sim$neretin

wwwth.itep.ru/$\sim$neretin
}
\end{document}